\documentclass[11pt]{article}
\usepackage{epsfig}

\def\be{\begin{equation}}
\def\ee{\end{equation}}
\def\bea{\begin{eqnarray}}
\def\eea{\end{eqnarray}}
\def\bes{\begin{eqnarray*}}
\def\ees{\end{eqnarray*}}

\def\nn{\nonumber}
\def\<{\langle}
\def\>{\rangle}
\def\lb{\label}
\def\bs{\setminus}

\def\R{{\bf R}}
\def\C{{\bf C}}
\def\Z{{\bf Z}}
\def\N{{\bf N}}
\def\U{{\bf U}}
\def\Q{{\bf Q}}

\def\ga{{\gamma}}

\def\th{{\theta}}

\def\Om{{\Omega}}

\def\dl{{\delta}}

\def\Sg{{\Sigma}}
\def\vf{{\varphi}}

\def\H{{\cal H}}

\def\J{{\cal J}}

\def\Nn{{\cal N}}

\def\Sp{{\rm Sp}}

\def\dm{{\rm \diamond}}

\def\td#1{\tilde{#1}}
\def\hb{\vrule height0.18cm width0.14cm $\,$}

\def\td#1{\tilde{#1}}
\def\mapright#1{\smash{\mathop{\longrightarrow}\limits^{#1}}}

\title{Stability of closed characteristics on compact convex\\ hypersurfaces in $\R^6$}
\author{Wei Wang\thanks{Partially supported by NNSF, RFDP of MOE of
China. E-mail: alexanderweiwang@yahoo.com.cn }\\
School of Mathematical Science \\ Peking University, Beijing 100871 \\
PEOPLES REPUBLIC OF CHINA \\ }
\date{}

\date{}

\begin{document}

\maketitle

\begin{abstract}
{\it In this paper, let $\Sigma\subset\R^{6}$ be a compact convex
hypersurface. We prove that if $\Sigma$ carries only finitely many
geometrically distinct closed characteristics, then at least two
of them must possess irrational mean indices. Moreover, if $\Sg$
carries exactly three geometrically distinct closed characteristics,
then at least two of them must be elliptic. }
\end{abstract}



{\bf Running title}: Stability of closed characteristics

\renewcommand{\theequation}{\thesection.\arabic{equation}}
\renewcommand{\thefigure}{\thesection.\arabic{figure}}

\setcounter{equation}{0}
\section{Introduction and main results}

In this paper, let $\Sigma$ be a fixed $C^3$ compact convex hypersurface
in $\R^{2n}$, i.e., $\Sigma$ is the boundary of a compact and strictly
convex region $U$ in $\R^{2n}$. We denote the set of all such hypersurfaces
by $\H(2n)$. Without loss of generality, we suppose $U$ contains the origin.
We consider closed characteristics $(\tau,y)$ on $\Sigma$, which are
solutions of the following problem
\be \left\{\matrix{\dot{y}=JN_{\Sigma}(y), \cr
               y(\tau)=y(0), \cr }\right. \lb{1.1}\ee
where $J=\left(\matrix{0 &-I_n\cr
                I_n  & 0\cr}\right)$,
$I_n$ is the identity matrix in $\R^n$, $\tau>0$, $N_\Sigma(y)$ is
the outward normal vector of $\Sigma$ at $y$ normalized by the
condition $N_{\Sigma}(y)\cdot y=1$. Here $a\cdot b$ denotes the
standard inner product of $a, b\in\R^{2n}$. A closed characteristic
$(\tau, y)$ is {\it prime}, if $\tau$ is the minimal period of $y$.
Two closed characteristics $(\tau, y)$ and $(\sigma, z)$ are {\it
geometrically distinct},  if $y(\R)\not= z(\R)$. We denote by
$\J(\Sg)$ and $\widetilde{\J}(\Sg)$ the set of all closed
characteristics $(\tau,\, y)$ on $\Sg$ with $\tau$ being the minimal
period of $y$ and the set of all geometrically distinct ones
respectively. Note that
$\J(\Sg)=\{\theta\cdot y\,|\, \theta\in S^1,\;y\; is\; prime\}$,
while $\widetilde{\J}(\Sg)=\J(\Sg)/S^1$, where the natural $S^1$-action
is defined by $\theta\cdot y(t)=y(t+\tau\theta),\;\;\forall \theta\in S^1,\,t\in\R$.

Let $j: \R^{2n}\rightarrow\R$ be the gauge function of $\Sigma$, i.e.,
$j(\lambda x)=\lambda$ for $x\in\Sigma$ and $\lambda\ge0$, then
$j\in C^3(\R^{2n}\setminus\{0\}, \R)\cap C^0(\R^{2n}, \R)$
and $\Sigma=j^{-1}(1)$. Fix a constant $\alpha\in(1,\,2)$ and
define the Hamiltonian function
$H_\alpha :\R^{2n}\rightarrow [0,\,+\infty)$ by
\be H_\alpha(x)=j(x)^\alpha,\qquad \forall x\in\R^{2n}.\lb{1.2}\ee
Then
$H_\alpha\in C^3(\R^{2n}\setminus\{0\}, \R)\cap C^1(\R^{2n}, \R)$
is convex and $\Sigma=H_\alpha^{-1}(1)$.
It is well known that the problem (\ref{1.1}) is equivalent to
the following given energy problem of the Hamiltonian system
\be
\left\{\matrix{\dot{y}(t)=JH_\alpha^\prime(y(t)),
             &&\quad H_\alpha(y(t))=1,\qquad \forall t\in\R. \cr
     y(\tau)=y(0). && \cr }\right. \lb{1.3}\ee
Denote by $\mathcal{J}(\Sigma, \,\alpha)$ the set of all solutions
$(\tau,\, y)$ of (\ref{1.3}) where $\tau$ is the minimal period of
$y$ and by $\widetilde{\mathcal{J}}(\Sigma, \,\alpha)$ the set of
all geometrically distinct solutions  of (\ref{1.3}). As above,
$\widetilde{\mathcal{J}}(\Sigma, \,\alpha)$ is obtained from
$\mathcal{J}(\Sigma, \,\alpha)$ by dividing the natural
$S^1$-action. Note that elements in $\mathcal{J}(\Sigma)$ and
$\mathcal{J}(\Sigma, \,\alpha)$ are one to one correspondent to each
other, similarly for $\widetilde{\J}(\Sg)$ and
$\widetilde{\mathcal{J}}(\Sigma, \,\alpha)$.

Let $(\tau,\, y)\in\mathcal{J}(\Sigma, \,\alpha)$. The fundamental
solution $\gamma_y : [0,\,\tau]\rightarrow \Sp(2n)$ with $\gamma_y(0)=I_{2n}$
of the linearized Hamiltonian system
\be \dot w(t)=JH_\alpha^{\prime\prime}(y(t))w(t),\qquad \forall t\in\R,\lb{1.4}\ee
is called the {\it associate symplectic path} of $(\tau,\, y)$.
The eigenvalues of $\gamma_y(\tau)$ are called {\it Floquet multipliers}
of $(\tau,\, y)$. By Proposition 1.6.13 of \cite{Eke3}, the Floquet multipliers
with their multiplicities of $(\tau,\, y)\in\mathcal{J}(\Sigma)$ do not depend on
the particular choice of the Hamiltonian function in (\ref{1.3}).
For any $M\in \Sp(2n)$, we define the {\it elliptic height } $e(M)$ of
$M$ to be the total algebraic multiplicity of all eigenvalues of $M$ on the
unit circle $\U=\{z\in\C|\; |z|=1\}$ in the complex plane $\C$.
Since $M$ is symplectic, $e(M)$ is even and $0\le e(M)\le 2n$.
As usual a $(\tau,\, y)\in\J(\Sg)$ is {\it elliptic}, if
$e(\gamma_y(\tau))=2n$. It is {\it non-degenerate}, if $1$ is a double
Floquet multiplier of it. It is  {\it hyperbolic}, if $1$ is a
double Floquet multiplier of it and $e(\gamma_y(\tau))=2$.
It is well known that these concepts are independent of the choice of $\alpha>1$.

For the existence and multiplicity of geometrically distinct closed
characteristics on convex compact hypersurfaces in $\R^{2n}$ we refer to
\cite{Rab1}, \cite{Wei1}, \cite{EkL1}, \cite{EkH1}, \cite{Szu1}, \cite{HWZ1},
\cite{LoZ1}, \cite{LLZ1}, and references therein. Note that recently in
\cite{WHL}, Wang, Hu and Long proved $\,^{\#}\td{\J}(\Sg)\ge 3$ for
every $\Sg\in\H(6)$.

On the stability problem, in \cite{Eke2} of Ekeland in 1986 and \cite{Lon2}
of Long in 1998, for any $\Sg\in\H(2n)$ the existence of at least one
non-hyperbolic closed characteristic on $\Sg$ was proved provided
$^\#\td{\J}(\Sg)<+\infty$. Ekeland proved also in \cite{Eke2} the existence
of at least one elliptic closed characteristic on $\Sg$ provided $\Sg\in\H(2n)$
is $\sqrt{2}$-pinched. In \cite{DDE1} of 1992, Dell'Antonio, D'Onofrio and
Ekeland proved the existence of at least one elliptic closed characteristic
on $\Sg$ provided $\Sg\in\H(2n)$ satisfies $\Sg=-\Sg$. In \cite{Lon3} of 2000,
Long proved that $\Sg\in\H(4)$ and $\,^{\#}\td{\J}(\Sg)=2$ imply that both of
the closed characteristics must be elliptic. In \cite{LoZ1} of 2002, Long and
Zhu further proved when $^\#\td{\J}(\Sg)<+\infty$, there exists at least one
elliptic closed characteristic and there are at least $[\frac{n}{2}]$ geometrically
distinct closed characteristics on $\Sg$ possessing irrational mean indices,
which are then non-hyperbolic. In the recent paper \cite{LoW1}, Long and Wang
proved that there exist at least two non-hyperbolic closed characteristic on
$\Sg\in\H(6)$ when $^\#\td{\J}(\Sg)<+\infty$. Motivated by these results, we prove
the following results in this paper:

{\bf Theorem 1.1.} {\it On every $\Sg\in\H(6)$ satisfying
$^\#\td{\J}(\Sg)<+\infty$, there exist at least two geometrically distinct
closed characteristics possessing irrational mean indices. }

{\bf Theorem 1.2.} {\it Suppose $^\#\td{\J}(\Sigma)=3$ for some
$\Sigma\in\H(6)$. Then there exist at least two elliptic closed
characteristics in $\td{\J}(\Sigma)$.}

The proofs of Theorems 1.1 and 1.2 are given in Section 3. Mainly ingredients in
the proofs inculde: the mean index identity for closed characteristics established
in \cite{WHL} recently, Morse inequality and the index iteration theory developed
by Long and his coworkers, specially the common index jump theorem of Long and Zhu
(Theorem 4.3 of \cite{LoZ1}, cf. Theorem 11.2.1 of \cite{Lon4}). In Section 2, we
review briefly the equivariant Morse theory and the mean index identity for closed
characteristics on compact convex hypersurfaces in $\R^{2n}$ developed in the
recent \cite{WHL}.

In this paper, let $\N$, $\N_0$, $\Z$, $\Q$, $\R$, and $\R^+$ denote
the sets of natural integers, non-negative integers, integers, rational
numbers, real numbers, and positive real numbers respectively.
Denote by $a\cdot b$ and $|a|$ the standard inner product and norm in
$\R^{2n}$. Denote by $\langle\cdot,\cdot\rangle$ and $\|\cdot\|$
the standard $L^2$-inner product and $L^2$-norm. For an $S^1$-space $X$, we denote
by $X_{S^1}$ the homotopy quotient of $X$ module the $S^1$-action, i.e.,
$X_{S^1}=S^\infty\times_{S^1}X$. We define the functions
\be \left\{\matrix{[a]=\max\{k\in\Z\,|\,k\le a\}, &
E(a)=\min\{k\in\Z\,|\,k\ge a\} , \cr
                   \varphi(a)=E(a)-[a],   \cr}\right. \lb{1.5}\ee
Specially, $\varphi(a)=0$ if $ a\in\Z\,$, and $\varphi(a)=1$ if $a\notin\Z\,$.
In this paper we use only $\Q$-coefficients for all homological modules.
For a $\Z_m$-space pair $(A, B)$, let
$H_{\ast}(A, B)^{\pm\Z_m}= \{\sigma\in H_{\ast}(A, B)\,|\,L_{\ast}\sigma=\pm \sigma\}$,
where $L$ is a generator of the $\Z_m$-action.

\setcounter{equation}{0}
\section{ Equivariant Morse theory for closed characteristics}

In the rest of this paper, we fix a $\Sg\in\H(2n)$ and assume the following
condition on $\Sg$:

\noindent (F) {\bf There exist only finitely many geometrically distinct
closed characteristics \\$\quad \{(\tau_j, y_j)\}_{1\le j\le k}$ on $\Sigma$. }

In this section, we review briefly the equivariant Morse theory
for closed characteristics on $\Sg$ developed in \cite{WHL} which will be
needed in Section 3 of this paper. All the details of proofs can be found
in \cite{WHL}.

Let $\hat{\tau}=\inf\{\tau_j|\;1\le j\le k\}$. Note that here
$\tau_j$'s are prime periods of $y_j$'s for $1\le j\le k$.
Then by \S2 of \cite{WHL}, for any $a>\hat{\tau}$, we can construct a function $\varphi_a\in C^\infty(\R,\R^+)$
which has $0$ as its unique critical point in $[0,\,+\infty)$ such that $\varphi_a$ is
strictly convex for $t\ge 0$. Moreover, $\frac{\varphi_a^\prime(t)}{t}$ is strictly
decreasing for $t> 0$ together with
$\lim_{t\rightarrow 0^+}\frac{\varphi_a^\prime(t)}{t}=1$
and $\varphi_a(0)=0=\varphi_a^\prime(0)$.
More precisely, we define $\varphi_a$ via Propositions 2.2
and 2.4 in \cite{WHL}. The precise dependence of
$\varphi_a$ on $a$ is explained in Remark 2.3 of \cite{WHL}.

Define the Hamiltonian
function $H_a(x)=a\varphi_a(j(x))$ and consider the fixed period problem
\be
\left\{\matrix{\dot{x}(t)=JH_a^\prime(x(t)), \cr
     x(1)=x(0).         \cr }\right. \lb{2.1}\ee
Then  $H_a\in C^3(\R^{2n}\setminus\{0\}, \R)\cap C^1(\R^{2n}, \R)$
is strictly convex. Solutions of (\ref{2.1}) are $x\equiv0$ and
$x=\rho y(\tau t)$ with
$\frac{\varphi_a^\prime(\rho)}{\rho}=\frac{\tau}{a}$,
where $(\tau, y)$ is a solution of (\ref{1.1}).
In particular, nonzero solutions of (\ref{2.1})
are one to one correspondent to solutions of (\ref{1.1})
with period $\tau<a$.

In the following, we use the Clarke-Ekeland dual action principle.
As usual, let $G_a$ be the Fenchel transform of $H_a$ defined by
$G_a(y)=\sup\{x\cdot y-H_a(x)\;|\; x\in \R^{2n}\}$. Then
$G_a\in C^2(\R^{2n}\bs\{0\},\R)\cap C^1(\R^{2n},\R)$ is strictly convex.
Let
\be L_0^2(S^1, \;\R^{2n})= \left\{u\in L^2([0, 1],\;\R^{2n})
   \left|\frac{}{}\right.\int_0^1u(t)dt=0\right\}.  \lb{2.2}\ee
Define a linear operator $M: L_0^2(S^1,\R^{2n})\to L_0^2(S^1,\R^{2n})$
by $\frac{d}{dt}Mu(t)=u(t)$, $\int_0^1Mu(t)dt=0$. The dual action
functional on $L_0^2(S^1, \;\R^{2n})$ is defined by
\be \Psi_a(u)=\int_0^1\left(\frac{1}{2}Ju\cdot Mu+G_a(-Ju)\right)dt.
   \lb{2.3}\ee
Then the functional $\Psi_a\in C^{1, 1}(L_0^2(S^1,\; \R^{2n}),\;\R)$
is bounded from below and satisfies the Palais-Smale condition. Suppose
$x$ is a solution of (\ref{2.1}). Then $u=\dot{x}$ is a critical point
of $\Psi_a$. Conversely, suppose $u$ is a critical point of $\Psi_a$.
Then there exists a unique $\xi\in\R^{2n}$ such that $Mu-\xi$ is a
solution of (\ref{2.1}). In particular, solutions of (\ref{2.1}) are in
one to one correspondence with critical points of $\Psi_a$. Moreover,
$\Psi_a(u)<0$ for every critical point $u\not= 0$ of $\Psi_a$.

Suppose $u$ is a nonzero critical point of $\Psi_a$. Then following
\cite{Eke3} the formal Hessian of $\Psi_a$ at $u$ is defined by
$$ Q_a(v,\; v)=\int_0^1 (Jv\cdot Mv+G_a^{\prime\prime}(-Ju)Jv\cdot Jv)dt, $$
which defines an orthogonal splitting $L_0^2=E_-\oplus E_0\oplus E_+$ of
$L_0^2(S^1,\; \R^{2n})$ into negative, zero and positive subspaces. The
index of $u$ is defined by $i(u)=\dim E_-$ and the nullity of $u$ is
defined by $\nu(u)=\dim E_0$. Let $u=\dot{x}$ be the critical point
of $\Psi_a$ such that $x$ corresponds to the closed characteristic
$(\tau,\,y)$ on $\Sigma$. Then the index $i(u)$ and the nullity
$\nu(u)$ defined above coincide with the Ekeland indices defined by I.
Ekeland in \cite{Eke1} and \cite{Eke3}. Specially $1\le \nu(u)\le 2n-1$
always holds.

We have a natural $S^1$-action on $L_0^2(S^1,\; \R^{2n})$ defined by
$\th\cdot u(t)=u(\th+t)$ for all $\th\in S^1$ and $t\in\R$. Clearly
$\Psi_a$ is $S^1$-invariant. For any $\kappa\in\R$, we denote by
\be \Lambda_a^\kappa=\{u\in L_0^2(S^1,\; \R^{2n})\;|\;\Psi_a(u)\le\kappa\}.
          \lb{2.4}\ee
For a critical point $u$ of $\Psi_a$, we denote by
\be \Lambda_a(u)=\Lambda_a^{\Psi_a(u)}
  =\{w\in L_0^2(S^1,\; \R^{2n}) \;|\; \Psi_a(w)\le\Psi_a(u)\}.\lb{2.5}\ee
Clearly, both sets are $S^1$-invariant. Since the
$S^1$-action preserves $\Psi_a$, if $u$ is a critical
point of $\Psi_a$, then the whole orbit $S^1\cdot u$ is formed by
critical points of $\Psi_a$. Denote by $crit(\Psi_a)$ the set of
critical points of $\Psi_a$. Note that by the condition (F),
the number of critical orbits of $\Psi_a$ is finite.
Hence as usual we can make the following definition.

{\bf Definition 2.1.} {\it Suppose $u$ is a nonzero critical
point of $\Psi_a$ and $\Nn$ is an $S^1$-invariant
open neighborhood of $S^1\cdot u$ such that
$crit(\Psi_a)\cap(\Lambda_a(u)\cap \Nn)=S^1\cdot u$. Then
the $S^1$-critical modules of $S^1\cdot u$ are defined by}
$$ C_{S^1,\; q}(\Psi_a, \;S^1\cdot u)
=H_{q}((\Lambda_a(u)\cap\Nn)_{S^1},\;
((\Lambda_a(u)\setminus S^1\cdot u)\cap\Nn)_{S^1}). $$

We have the following proposition for critical modules.

{\bf Proposition 2.2.} (Proposition 3.2 of \cite{WHL}) {\it The critical
module $C_{S^1,\;q}(\Psi_a, \;S^1\cdot u)$ is independent of $a$ in the
sense that if $x_i$ are solutions of (\ref{2.1}) with Hamiltonian
functions $H_{a_i}(x)\equiv a_i\varphi_{a_i}(j(x))$ for $i=1$ and $2$
respectively such that both $x_1$ and $x_2$ correspond to the same closed
characteristic $(\tau, y)$ on $\Sigma$. Then we have}
$$ C_{S^1,\; q}(\Psi_{a_1}, \;S^1\cdot\dot {x}_1) \cong
  C_{S^1,\; q}(\Psi_{a_2}, \;S^1\cdot \dot {x}_2), \quad \forall q\in \Z. $$

Now let $u\neq 0$ be a critical point of $\Psi_a$ with multiplicity
$mul(u)=m$, i.e., $u$ corresponds to a closed characteristic
$(m\tau, y)\subset\Sigma$ with $(\tau, y)$ being prime. Hence
$u(t+\frac{1}{m})=u(t)$ holds for all $t\in \R$ and the orbit of $u$,
namely, $S^1\cdot u\cong S^1/\Z_m\cong S^1$.
Let $f: N(S^1\cdot u)\rightarrow S^1\cdot u$ be the normal
bundle of $S^1\cdot u$ in $L_0^2(S^1,\; \R^{2n})$ and let
$f^{-1}(\theta\cdot u)=N(\theta\cdot u)$ be the fibre over
$\theta\cdot u$, where $\theta\in S^1$.
Let $DN(S^1\cdot u)$ be the $\varrho$-disk bundle of $N(S^1\cdot u)$
for some $\varrho>0$ sufficiently small, i.e.,
$DN(S^1\cdot u)=\{\xi\in N(S^1\cdot u)\;| \; \|\xi\|<\varrho\}$
and let $DN(\theta\cdot u)=f^{-1}(\th\cdot u)\cap DN(S^1\cdot u)$
be the disk over $\theta\cdot u$. Clearly, $DN(\theta\cdot u)$ is
$\Z_m$-invariant and we have $DN(S^1\cdot u)=DN(u)\times_{\Z_m}S^1$,
where the $Z_m$-action is given by
$$ (\th, v, t)\in \Z_m\times DN(u)\times S^1\mapsto
        (\th\cdot v, \;\theta^{-1}t)\in DN(u)\times S^1. $$
Hence for an $S^1$-invariant subset $\Gamma$ of $DN(S^1\cdot u)$,
we have $\Gamma/S^1=(\Gamma_u\times_{\Z_m}S^1)/S^1=\Gamma_u/\Z_m$,
where $\Gamma_u=\Gamma\cap DN(u)$. Since $\Psi_a$ is not $C^2$
on $L_0^2(S^1,\; \R^{2n})$, we need to use a finite dimensional
approximation introduced by Ekeland in order to apply Morse theory.
More precisely, we can construct a finite dimensional submanifold
$\Gamma(\iota)$ of $L_0^2(S^1,\; \R^{2n})$ which admits a $\Z_\iota$-action
with $m|\iota$. Moreover $\Psi_a$ and $\Psi_a|_{\Gamma(\iota)}$ have
the same critical points. $\Psi_a|_{\Gamma(\iota)}$ is $C^2$ in a small
tubular neighborhood of the critical orbit $S^1\cdot u$ and the
Morse index and nullity of its critical points coincide with those of
the corresponding critical points of $\Psi_a$.  Let
\be D_\iota N(S^1\cdot u)=DN(S^1\cdot u)\cap\Gamma(\iota), \quad
D_\iota N(\theta\cdot u)=DN(\theta\cdot u)\cap\Gamma(\iota). \lb{2.6}\ee
Then we have
\be C_{S^1,\; \ast}(\Psi_a, \;S^1\cdot u)
\cong H_\ast(\Lambda_a(u)\cap D_\iota N(u),\;
    (\Lambda_a(u)\setminus\{u\})\cap D_\iota N(u))^{\Z_m}. \lb{2.7}\ee
Now we can apply the results of Gromoll and Meyer in \cite{GrM1}
to the manifold $D_{p\iota}N(u^p)$ with $u^p$ as its unique critical point,
where $p\in\N$. Then $mul(u^p)=pm$ is the multiplicity of $u^p$ and the
isotropy group $\Z_{pm}\subseteq S^1$ of $u^p$ acts on $D_{p\iota}N(u^p)$
by isometries. According to Lemma 1 of \cite{GrM1}, we have a
$\Z_{pm}$-invariant decomposition of $T_{u^p}(D_{p\iota}N(u^p))$
$$ T_{u^p}(D_{p\iota}N(u^p))
=V^+\oplus V^-\oplus V^0=\{(x_+, x_-, x_0)\}  $$
with $\dim V^-=i(u^p)$, $\dim V^0=\nu(u^p)-1$ and a
$\Z_{pm}$-invariant neighborhood $B=B_+\times B_-\times B_0$
for $0$ in $T_{u^p}(D_{p\iota}N(u^p))$ together with two $Z_{pm}$-invariant
diffeomorphisms
$$\Phi :B=B_+\times B_-\times B_0\rightarrow
\Phi(B_+\times B_-\times B_0)\subset D_{p\iota}N(u^p)$$
and
$$ \eta : B_0\rightarrow W(u^p)\equiv\eta(B_0)\subset D_{p\iota}N(u^p)$$
such that $\Phi(0)=\eta(0)=u^p$ and
\be \Psi_a\circ\Phi(x_+,x_-,x_0)=|x_+|^2 - |x_-|^2 + \Psi_a\circ\eta(x_0),
    \lb{2.8}\ee
with $d(\Psi_a\circ \eta)(0)=d^2(\Psi_a\circ\eta)(0)=0$.
As \cite{GrM1}, we call $W(u^p)$ a local {\it characteristic manifold} and
$U(u^p)=B_-$ a local {\it negative disk} at $u^p$.
By the proof of Lemma 1 of \cite{GrM1},
$W(u^p)$ and $U(u^p)$ are $\Z_{pm}$-invariant. Then we have
\bea
&& H_\ast(\Lambda_a(u^p)\cap D_{p\iota}N(u^p),\;
  (\Lambda_a(u^p)\setminus\{u^p\})\cap D_{p\iota}N(u^p)) \nn\\
&&\qquad = H_\ast (U(u^p),\;U(u^p)\setminus\{u^p\}) \otimes
H_\ast(W(u^p)\cap \Lambda_a(u^p),\; (W(u^p)\setminus\{u^p\})\cap \Lambda_a(u^p)),
  \lb{2.9}\eea
where \be H_q(U(u^p),U(u^p)\setminus\{u^p\} )
    = \left\{\matrix{\Q, & {\rm if\;}q=i(u^p),  \cr
                      0, & {\rm otherwise}. \cr}\right.  \lb{2.10}\ee
Now we have the following proposition.

{\bf Proposition 2.3.} (Proposition 3.10 of \cite{WHL}) {\it Let $u\neq 0$ be
a critical point of $\Psi_a$ with $mul(u)=1$. Then for all $p\in\N$
and $q\in\Z$, we have
\be C_{S^1,\; q}(\Psi_a, \;S^1\cdot u^p)\cong
\left(\frac{}{}H_{q-i(u^p)}(W(u^p)\cap \Lambda_a(u^p),\;
(W(u^p)\setminus\{u^p\})\cap \Lambda_a(u^p))\right)^{\beta(u^p)\Z_p},
  \lb{2.11}\ee
where $\beta(u^p)=(-1)^{i(u^p)-i(u)}$. Thus
\be C_{S^1,\; q}(\Psi_a, \;S^1\cdot u^p)=0, \quad {\rm for}\;\;
   q<i(u^p) \;\;{\rm or}\;\;q>i(u^p)+\nu(u^p)-1. \lb{2.12}\ee
In particular, if $u^p$ is
non-degenerate, i.e., $\nu(u^p)=1$, then}
\be C_{S^1,\; q}(\Psi_a, \;S^1\cdot u^p)
    = \left\{\matrix{\Q, & {\rm if\;}q=i(u^p)\;{\rm and\;}\beta(u^p)=1,  \cr
                      0, & {\rm otherwise}. \cr}\right.  \lb{2.13}\ee

We make the following definition

{\bf Definition 2.4.} {\it Let $u\neq 0$ be a critical
point of $\Psi_a$ with $mul(u)=1$.
Then for all $p\in\N$ and $l\in\Z$, let
\bea
k_{l, \pm 1}(u^p)&=&\dim\left(\frac{}{}H_l(W(u^p)\cap \Lambda_a(u^p),\;
(W(u^p)\setminus\{u^p\})\cap \Lambda_a(u^p))\right)^{\pm\Z_p}, \nn\\
k_l(u^p)&=&\dim\left(\frac{}{}H_l(W(u^p)\cap \Lambda_a(u^p),
(W(u^p)\setminus\{u^p\})\cap \Lambda_a(u^p))\right)^{\beta(u^p)\Z_p}. \nn\eea
$k_l(u^p)$'s are called critical type numbers of $u^p$. }

We have the following properties for critical type numbers

{\bf Proposition 2.5.} (Proposition 3.13 of \cite{WHL}) {\it Let $u\neq 0$
be a critical point of $\Psi_a$ with $mul(u)=1$. Then there exists a minimal
$K(u)\in \N$ such that
$$ \nu(u^{p+K(u)})=\nu(u^p),\quad i(u^{p+K(u)})-i(u^p)\in 2\Z, $$
and $k_l(u^{p+K(u)})=k_l(u^p)$ for all $p\in \N$ and $l\in\Z$. We call
$K(u)$ the minimal period of critical modules of iterations of the
functional $\Psi_a$ at $u$.   }

For a closed characteristic $(\tau,y)$ on $\Sigma$, we denote by
$y^m\equiv (m\tau, y)$ the $m$-th iteration of $y$ for $m\in\N$.
Let $a>\tau$ and choose $\vf_a$ as above. Determine $\rho$ uniquely by
$\frac{\vf_a'(\rho)}{\rho}=\frac{\tau}{a}$. Let $x=\rho y(\tau t)$ and
$u=\dot{x}$. Then we define the index $i(y^m)$ and nullity $\nu(y^m)$
of $(m\tau,y)$ for $m\in\N$ by
$$ i(y^m)=i(u^m), \qquad \nu(y^m)=\nu(u^m). $$
These indices are independent of $a$ when $a$ tends to infinity.
Now the mean index of $(\tau,y)$ is defined by
$$ \hat{i}(y)=\lim_{m\rightarrow\infty}\frac{i(y^m)}{m}. $$
Note that $\hat{i}(y)>2$ always holds which was proved by Ekeland and
Hofer in \cite{EkH1} of 1987 (cf. Corollary 8.3.2 and Lemma 15.3.2
of \cite{Lon4} for a different proof).

By Proposition 2.2, we can define the critical type numbers $k_l(y^m)$
of $y^m$ to be $k_l(u^m)$, where $u^m$ is the critical point of $\Psi_a$
corresponding to $y^m$. We also define $K(y)=K(u)$. Then we have

{\bf Proposition 2.6.} {\it We have $k_l(y^m)=0$ for $l\notin [0, \nu(y^m)-1]$
and it can take only values $0$ or $1$ when $l=0$ or $l=\nu(y^m)-1$.
Moreover, the following properties hold (cf. Lemma 3.10 of \cite{BaL1},
\cite{Cha1} and \cite{MaW1}):

(i) $k_0(y^m)=1$ implies $k_l(y^m)=0$ for $1\le l\le \nu(y^m)-1$.

(ii) $k_{\nu(y^m)-1}(y^m)=1$ implies $k_l(y^m)=0$ for $0\le l\le \nu(y^m)-2$.

(iii) $k_l(y^m)\ge 1$ for some $1\le l\le \nu(y^m)-2$ implies
$k_0(y^m)=k_{\nu(y^m)-1}(y^m)=0$.

(iv) If $\nu(y^m)\le 3$, then at most one of the $k_l(y^m)$'s for
$0\le l\le \nu(y^m)-1$ can be non-zero.

(v) If $i(y^m)-i(y)\in 2\Z+1$ for some $m\in\N$, then $k_0(y^m)=0$.}

{\bf Proof.} By Definition 2.4 we have
$$ k_l(y^m)\le \dim H_l(W(u^m)\cap \Lambda_a(u^m),\;
(W(u^m)\setminus\{u^m\})\cap \Lambda_a(u^m))\equiv \eta_l(y^m). $$
Then from Corollary 1.5.1 of \cite{Cha1} or Corollary 8.4 of \cite{MaW1},
(i)-(iv) hold.

For (v), if $\eta_0(y^m)=0$, then (v) follows directly from Definition 2.4.

By Corollary 8.4 of \cite{MaW1}, $\eta_0(y^m)=1$ if and only if $u^m$ is a
local minimum in the local characteristic manifold $W(u^m)$. Hence
$(W(u^m)\cap \Lambda_a(u^m),\;(W(u^m)\setminus\{u^m\})\cap
\Lambda_a(u^m))=(\{u^m\},\; \emptyset)$. By Definition 2.4, we have:
\bea k_{0, +1}(u^m)
&=& \dim H_0(W(u^m)\cap \Lambda_a(u^m),\;
    (W(u^m)\setminus\{u^m\})\cap \Lambda_a(u^m))^{+\Z_m}\nn\\
&=& \dim H_0(\{u^m\})^{+\Z_m}\nn\\
&=& 1.   \nn\eea
This implies $k_0(u^m)=k_{0, -1}(u^m)=0$. \hfill\hb

For a closed characteristic $(\tau, y)$ on $\Sigma$, we define as in \cite{WHL}
\be \hat\chi(y)=\frac{1}{K(y)}
  \sum_{1\le m\le K(y)\atop 0\le l\le 2n-2}
  (-1)^{i(y^{m})+l}k_l(y^{m}). \lb{2.14}\ee
In particular, if all $y^m$'s are non-degenerate, then by Proposition 2.3 we have
\be \hat\chi(y)
    = \left\{\matrix{(-1)^{i(y)}, & {\rm if\;\;} i(y^2)-i(y)\in 2\Z,  \cr
           \frac{(-1)^{i(y)}}{2}, & {\rm otherwise}. \cr}\right.  \lb{2.15}\ee

We have the following mean index identity for closed characteristics.

{\bf Theorem 2.7.} (Theorem 1.2 of \cite{WHL}) {\it Suppose
$\Sigma\in \H(2n)$ satisfies $\,^{\#}\widetilde{\J}(\Sg)<+\infty$. Denote all
the geometrically distinct closed characteristics by
$\{(\tau_j,y_j)\}_{1\le j\le k}$. Then the following identity holds }
$$ \sum_{1\le j\le k}\frac{\hat{\chi}(y_j)}{\hat{i}(y_j)}=\frac{1}{2}. $$

Let $\Psi_a$ be the functional defined by (\ref{2.3}) for
some $a\in\R$ large enough and let $\varepsilon>0$ be small
enough such that $[-\varepsilon, +\infty)\setminus\{0\}$ contains no critical
values of $\Psi_a$. Denote by $I_a$ the greatest integer in $\N_0$ such that
$I_a<i(\tau, y)$ hold for all closed characteristics $(\tau,\, y)$ on $\Sigma$ with
$\tau\ge a$. Then by Section 5 of \cite{WHL}, we have
\be H_{S^1,\; q}(\Lambda_a^{-\varepsilon} ) \cong H_{S^1,\; q}( \Lambda_a^\infty)
  \cong H_q(CP^\infty), \quad \forall q<I_a.  \lb{2.16}\ee
For any $q\in\Z$, let
\be  M_q(\Lambda_a^{-\varepsilon})
  =\sum_{1\le j\le k,\,1\le m_j<a/\tau_j} \dim C_{S^1,\;q}(\Psi_a, \;S^1\cdot u_j^{m_j}).
  \lb{2.17} \ee
Then the equivariant Morse inequalities for the space $\Lambda_a^{-\varepsilon}$
yield
\bea M_q(\Lambda_a^{-\varepsilon})
       &\ge& b_q(\Lambda_a^{-\varepsilon}),\lb{2.18}\\
M_q(\Lambda_a^{-\varepsilon}) &-& M_{q-1}(\Lambda_a^{-\varepsilon})
    + \cdots +(-1)^{q}M_0(\Lambda_a^{-\varepsilon}) \nn\\
&\ge& b_q(\Lambda_a^{-\varepsilon}) - b_{q-1}(\Lambda_a^{-\varepsilon})
   + \cdots + (-1)^{q}b_0(\Lambda_a^{-\varepsilon}), \lb{2.19}\eea
where $b_q(\Lambda_a^{-\varepsilon})=\dim H_{S^1,\; q}(\Lambda_a^{-\varepsilon})$.
Now we have the following Morse inequalities for closed characteristics

{\bf Theorem 2.8.} {\it Let $\Sigma\in \H(2n)$ satisfy
$\,^{\#}\widetilde{\J}(\Sg)<+\infty$. Denote all the geometrically distinct closed
characteristics by $\{(\tau_j,\; y_j)\}_{1\le j\le k}$. Let
\bea
M_q&=&\lim_{a\rightarrow+\infty}M_q(\Lambda_a^{-\varepsilon}),\quad
                  \forall q\in\Z,\lb{2.20}\\
b_q &=& \lim_{a\rightarrow+\infty}b_q(\Lambda_a^{-\varepsilon})=
\left\{\matrix{1, & {\rm if\;}q\in 2\N_0,  \cr
                      0, & {\rm otherwise}. \cr}\right.  \lb{2.21}
\eea
Then we have}
\bea M_q &\ge& b_q,\lb{2.22}\\
 M_q-M_{q-1}+\cdots +(-1)^{q}M_0 &\ge& b_q-b_{q-1}+\cdots +(-1)^{q}b_0,
   \qquad\forall \;q\in\Z. \lb{2.23}\eea

{\bf Proof.} As we have mentioned before, $\hat i(y_j)>2$ holds for $1\le j\le k$.
Hence the Ekeland index satisfies $i(y_j^m)=i(u_j^m)\to\infty$ as $m\to\infty$ for
$1\le j\le k$. Note that $I_a\to +\infty$ as $a\to +\infty$. Now fix a $q\in\Z$ and
a sufficiently great $a>0$. By Propositions 2.2, 2.3 and (\ref{2.17}),
$M_i(\Lambda_a^{-\varepsilon})$ is invariant for all $a>A_q$ and $0\le i\le q$,
where $A_q>0$ is some constant. Hence (\ref{2.20}) is meaningful. Now for any $a$
such that $I_a>q$, (\ref{2.16})-(\ref{2.19}) imply that (\ref{2.21})-(\ref{2.23})
hold. \hfill\hb

\setcounter{equation}{0}
\section{Proofs of the main theorems }

In this section, we give proofs of Theorems 1.1 and 1.2 by using the mean index
identity of \cite{WHL}, Morse inequality and the index iteration theory developed
by Long and his coworkers.

As Definition 1.1 of \cite{LoZ1}, we define

{ \bf Definition  3.1.} For $\alpha\in(1,2)$, we define a map
$\varrho_n\colon\H(2n)\to\N\cup\{ +\infty\}$
\be \varrho_n(\Sg)
= \left\{\matrix{+\infty, & {\rm if\;\;}^\#\mathcal{V}(\Sigma,\alpha)=+\infty, \cr
\min\left\{[\frac{i(x,1) + 2S^+(x) - \nu(x,1)+n}{2}]\,
\left|\frac{}{}\right.\,(\tau,x)\in\mathcal{V}_\infty(\Sigma, \alpha)\right\},
 & {\rm if\;\;} ^\#\mathcal{V}(\Sigma, \alpha)<+\infty, \cr}\right.  \lb{3.1}\ee
where $\mathcal{V}(\Sigma,\alpha)$ and $\mathcal{V}_\infty(\Sigma,\alpha)$ are
variationally visible and infinite variationally visible sets respectively given
by Definition 1.4 of \cite{LoZ1} (cf. Definition 15.3.3 of \cite{Lon4}).

{\bf Theorem 3.2.} (cf. Theorem 15.1.1 of \cite{Lon4}) {\it Suppose
$(\tau,y)\in \J(\Sigma)$. Then we have
\be i(y^m)\equiv i(m\tau ,y)=i(y, m)-n,\quad \nu(y^m)\equiv\nu(m\tau, y)=\nu(y, m),
       \qquad \forall m\in\N, \lb{3.2}\ee
where $i(y, m)$ and $\nu(y, m)$ are the Maslov-type index and nullity of
$(m\tau, y)$ defined by Conley, Zehnder and Long (cf. \S 5.4 of \cite{Lon4}).}

Recall that for a principal $U(1)$-bundle $E\to B$, the Fadell-Rabinowitz index
(cf. \cite{FaR1}) of $E$ is defined to be $\sup\{k\;|\, c_1(E)^{k-1}\not= 0\}$,
where $c_1(E)\in H^2(B,\Q)$ is the first rational Chern class. For a $U(1)$-space,
i.e., a topological space $X$ with a $U(1)$-action, the Fadell-Rabinowitz index is
defined to be the index of the bundle $X\times S^{\infty}\to X\times_{U(1)}S^{\infty}$,
where $S^{\infty}\to CP^{\infty}$ is the universal $U(1)$-bundle.

As in P.199 of \cite{Eke3}, choose some $\alpha\in(1,\, 2)$ and associate with $U$
a convex function $H$ such that $H(\lambda x)=\lambda^\alpha H(x)$ for $\lambda\ge 0$.
Consider the fixed period problem
\be \left\{\matrix{\dot{x}(t)=JH^\prime(x(t)), \cr
     x(1)=x(0).         \cr }\right. \lb{3.3}\ee

Define
\be L_0^{\frac{\alpha}{\alpha-1}}(S^1,\R^{2n})
  =\{u\in L^{\frac{\alpha}{\alpha-1}}(S^1,\R^{2n})\,|\,\int_0^1udt=0\}. \lb{3.4}\ee
The corresponding Clarke-Ekeland dual action functional is defined by
\be \Phi(u)=\int_0^1\left(\frac{1}{2}Ju\cdot Mu+H^{\ast}(-Ju)\right)dt,
    \qquad \forall\;u\in L_0^{\frac{\alpha}{\alpha-1}}(S^1,\R^{2n}), \lb{3.5}\ee
where $Mu$ is defined by $\frac{d}{dt}Mu(t)=u(t)$ and $\int_0^1Mu(t)dt=0$,
$H^\ast$ is the Fenchel transform of $H$ defined in \S2.

For any $\kappa\in\R$, we denote by
\be \Phi^{\kappa-}=\{u\in L_0^{\frac{\alpha}{\alpha-1}}(S^1,\R^{2n})\;|\;
             \Phi(u)<\kappa\}. \lb{3.6}\ee
Then as in P.218 of \cite{Eke3}, we define
\be c_i=\inf\{\delta\in\R\;|\: \hat I(\Phi^{\delta-})\ge i\},\lb{3.7}\ee
where $\hat I$ is the Fadell-Rabinowitz index given above. Then by Proposition 3
in P.218 of \cite{Eke3}, we have

{\bf Proposition 3.3.} {\it Every $c_i$ is a critical value of $\Phi$. If
$c_i=c_j$ for some $i<j$, then there are infinitely many geometrically
distinct closed characteristics on $\Sg$.}

As in Definition 2.1, we define the following

{\bf Definition 3.4.} {\it Suppose $u$ is a nonzero critical
point of $\Phi$, and $\Nn$ is an $S^1$-invariant
open neighborhood of $S^1\cdot u$ such that
$crit(\Phi)\cap(\Lambda(u)\cap \Nn)=S^1\cdot u$. Then
the $S^1$-critical modules of $S^1\cdot u$ is defined by
\bea C_{S^1,\; q}(\Phi, \;S^1\cdot u)
=H_{q}((\Lambda(u)\cap\Nn)_{S^1},\;
((\Lambda(u)\setminus S^1\cdot u)\cap\Nn)_{S^1}),\lb{3.8}
\eea
where $\Lambda(u)=\{w\in L_0^{\frac{\alpha}{\alpha-1}}(S^1,\R^{2n})\;|\;
\Phi(w)\le\Phi(u)\}$.}

Comparing with Theorem 4 in P.219 of \cite{Eke3}, we have the following

{\bf Proposition 3.5.} {\it For every $i\in\N$, there exists a point
$u\in L_0^{\frac{\alpha}{\alpha-1}}(S^1,\R^{2n})$ such that}
\bea
&& \Phi^\prime(u)=0,\quad \Phi(u)=c_i, \lb{3.9}\\
&& C_{S^1,\; 2(i-1)}(\Phi, \;S^1\cdot u)\neq 0. \lb{3.10}\eea

{\bf Proof.} By Lemma 8 in P.206 of \cite{Eke3}, we can use
Theorem 1.4.2 of \cite{Cha1} in the equivariant form  to obtain
\be H_{S^1,\,\ast}(\Phi^{c_i+\epsilon},\;\Phi^{c_i-\epsilon})
=\bigoplus_{\Phi(u)=c_i}C_{S^1,\; \ast}(\Phi, \;S^1\cdot u),\lb{3.11}\ee
for $\epsilon$  small enough such that the interval
$(c_i-\epsilon,\,c_i+\epsilon)$ contains no critical values of $\Phi$
except $c_i$.

Similar to P.431 of \cite{EkH1}, we have
\be H^{2(i-1)}((\Phi^{c_i+\epsilon})_{S^1},\,(\Phi^{c_i-\epsilon})_{S^1})
\mapright{q^\ast} H^{2(i-1)}((\Phi^{c_i+\epsilon})_{S^1} )
\mapright{p^\ast}H^{2(i-1)}((\Phi^{c_i-\epsilon})_{S^1}),
\lb{3.12}\ee
where $p$ and $q$ are natural inclusions. Denote by
$f: (\Phi^{c_i+\epsilon})_{S^1}\rightarrow CP^\infty$ a classifying map and let
$f^{\pm}=f|_{(\Phi^{c_i\pm\epsilon})_{S^1}}$. Then clearly each
$f^{\pm}: (\Phi^{c_i\pm\epsilon})_{S^1}\rightarrow CP^\infty$ is a classifying
map on $(\Phi^{c_i\pm\epsilon})_{S^1}$. Let $\eta \in H^2(CP^\infty)$ be the first
universal Chern class.

By definition of $c_i$, we have $\hat I(\Phi^{c_i-\epsilon})< i$, hence
$(f^-)^\ast(\eta^{i-1})=0$. Note that
$p^\ast(f^+)^\ast(\eta^{i-1})=(f^-)^\ast(\eta^{i-1})$.
Hence the exactness of (\ref{3.12}) yields a
$\sigma\in H^{2(i-1)}((\Phi^{c_i+\epsilon})_{S^1},\,(\Phi^{c_i-\epsilon})_{S^1})$
such that $q^\ast(\sigma)=(f^+)^\ast(\eta^{i-1})$.
Since $\hat I(\Phi^{c_i+\epsilon})\ge i$, we have $(f^+)^\ast(\eta^{i-1})\neq 0$.
Hence $\sigma\neq 0$, and then
$$H^{2(i-1)}_{S^1}(\Phi^{c_i+\epsilon},\Phi^{c_i-\epsilon})=
H^{2(i-1)}((\Phi^{c_i+\epsilon})_{S^1},\,(\Phi^{c_i-\epsilon})_{S^1})\neq 0. $$
Now the proposition follows from (\ref{3.11}) and the universal coefficient
theorem. \hfill\hb

{\bf Proposition 3.6.} {\it Suppose $u$ is the critical point of $\Phi$ found
in Proposition 3.5. Then we have
\be C_{S^1,\; 2(i-1)}(\Psi_a, \;S^1\cdot u_a)\neq 0, \lb{3.13}\ee
where $\Psi_a$ is given by (\ref{2.3}) and $u_a\in L_0^2(S^1,\;\R^{2n})$
is its critical point corresponding to $u$ in the natural sense.}

{\bf Proof.} Fix this $u$, we  modify the function $H$ only in a small
neighborhood $\Omega$ of $0$ as in \cite{Eke1} so that the corresponding
orbit of $u$ does not enter $\Omega$ and the resulted function $\widetilde{H}$
satisfies similar properties as Definition 1 in P. 26
of \cite{Eke1} by just replacing $\frac{3}{2}$ there by $\alpha$.
Define the dual action functional
$\widetilde{\Phi}:L_0^{\frac{\alpha}{\alpha-1}}(S^1,\R^{2n})\to\R$ by
\be \widetilde{\Phi}(v)=\int_0^1\left(\frac{1}{2}Jv\cdot
   Mv+\widetilde{H}^{\ast}(-Jv)\right)dt, \lb{3.14}\ee
since clearly $\Phi$ and $\widetilde{\Phi}$ are $C^1$ close to each other.
Then by the continuity of critical modules (cf. Theorem 8.8 of \cite{MaW1} or
Theorem 1.5.6 in P.53 of \cite{Cha1}, which can be easily generalized to the
equivariant sense) for the $u$ in the proposition, we have
\be C_{S^1,\; \ast}(\Phi, \;S^1\cdot u)\cong C_{S^1,\; \ast}(\widetilde{\Phi},
    \;S^1\cdot u).\lb{3.15}\ee

Using a finite dimensional approximation as in Lemma 3.9 of \cite{Eke1},
we have
\be C_{S^1,\; \ast}(\widetilde{\Phi}, \;S^1\cdot u)
\cong H_\ast(\widetilde{\Lambda}(u)\cap D_\iota N(u),\;
    (\widetilde{\Lambda}(u)\setminus\{u\})\cap D_\iota N(u))^{\Z_m}, \lb{3.16}\ee
where $\widetilde{\Lambda}(u)=\{w\in L_0^{\frac{\alpha}{\alpha-1}}(S^1,\R^{2n})\;|\;
\widetilde{\Phi}(w)\le\widetilde{\Phi}(u)\}$ and $D_\iota N(u)$ is a
$\Z_m$-invariant finite dimensional disk transversal to $S^1\cdot u$ at $u$
(cf. Lemma 3.9 of \cite{WHL}), $m$ is the multiplicity of $u$.

By Lemma 3.9 of \cite{WHL}, we have
\be C_{S^1,\; \ast}(\Psi_a, \;S^1\cdot u_a)
\cong H_\ast(\Lambda_a(u_a)\cap D_\iota N(u_a),\;
    (\Lambda_a(u_a)\setminus\{u_a\})\cap D_\iota N(u_a))^{\Z_m}.\lb{3.17}\ee
By the construction of $H_a$ in \cite{WHL}, $H_a=\widetilde{H}$ in a
$L^\infty$-neighborhood of $S^1\cdot u$. We remark here that multiplying $H$ by a constant
will not affect the corresponding critical modules, i.e., the corresponding
critical orbits have isomorphic critical modules. Hence we can assume
$H_a=H$ in a $L^\infty$-neighborhood of $S^1\cdot u$ and then the above
conclusion. Hence $\Psi_a$ and $\widetilde{\Phi}$ coincide
in a $L^\infty$-neighborhood of $S^1\cdot u$. Note also by Lemma 3.9 of
\cite{Eke1}, the two finite dimensional approximations are actually the same.
Hence we have
\bea
&& H_\ast(\widetilde{\Lambda}(u)\cap D_\iota N(u),\;
   (\widetilde{\Lambda}(u)\setminus\{u\})\cap D_\iota N(u))^{\Z_m}\nn\\
&&\quad\cong H_\ast(\Lambda_a(u_a)\cap D_\iota N(u_a),\;
    (\Lambda_a(u_a)\setminus\{u_a\})\cap D_\iota N(u_a))^{\Z_m}.\lb{3.18}\eea
Now the proposition follows from Proposition 3.5 and (\ref{3.16})-(\ref{3.18}).
\hfill\hb

Now we can give:

{\bf Proof of Theorem 1.1.} By the assumption (F) at the beginning of Section 2,
we denote by $\{(\tau_j, y_j)\}_{1\le j\le k}$ all the  geometrically distinct
closed characteristics on $\Sg$, and by $\ga_j\equiv \gamma_{y_j}$ the associated
symplectic path of $(\tau_j,\,y_j)$ on $\Sg$ for $1\le j\le k$. Then by Lemma
15.2.4 of \cite{Lon4}, there exist $P_j\in \Sp(6)$ and $M_j\in \Sp(4)$ such
that
\be \ga_j(\tau_j)=P_j^{-1}(N_1(1,\,1)\dm M_j)P_j, \quad\forall\; 1\le j\le k,
   \lb{3.19}\ee
where recall $N_1(1,b)=\left(\matrix{1 & b\cr
                         0 &  1\cr}\right)$ for $b\in\R$.

Without loss of generality, by Theorem 1.3 of \cite{LoZ1} (cf. Theorem 15.5.2 of
\cite{Lon4}), we may assume that $(\tau_1,y_1)$ has irrational mean index. Hence
by Theorem 8.3.1 and Corollary 8.3.2  of \cite{Lon4}, $M_1\in \Sp(4)$ in (\ref{3.19})
can be connected to $R(\th_1)\dm Q_1$ within $\Om^0(M_1)$ for some
$\frac{\th_1}{\pi}\notin\Q$ and $Q_1\in \Sp(2)$, where
$R(\th)=\left(\matrix{\cos\th & -\sin\th\cr
                      \sin\th &  \cos\th\cr}\right)$ for $\th\in\R$. Here we
use notations from Definition 1.8.5 and Theorem 1.8.10 of \cite{Lon4}. By
Theorem 2.7, the following identity holds
\be  \frac{\hat{\chi}(y_1)}{\hat{i}(y_1)} +
    \sum_{2\le j\le k}\frac{\hat{\chi}(y_j)}{\hat{i}(y_j)}=\frac{1}{2}.
\lb{3.20}\ee
Now we have the following four cases according to the classification
of basic norm forms (cf. Definition 1.8.9 of \cite{Lon4}).

{\bf Case 1.} {\it $Q_1=R(\th_2)$ with $\frac{\th_2}{\pi}\notin\Q$ or
$Q_1=D(\pm 2)\equiv\left(\matrix{ \pm 2 & 0\cr
                  0 & \pm\frac{1}{2}\cr}\right)$}.

\smallskip

In this case, by Theorems 8.1.6 and 8.1.7 of \cite{Lon4}, we have
$\nu(y_1^m)\equiv 1$, i.e., $y_1^m$ is non-degenerate for all $m\in\N$.
Hence it follows from (\ref{2.15}) that $\hat\chi(y_1)\neq 0$.
Now (\ref{3.20}) implies that at least one of the $y_j$'s for $2\le j\le k$
must have irrational mean index. Hence the theorem holds.

{\bf Case 2.} {\it $Q_1=N_1(1,b)$ with $b=\pm 1,\, 0$}.

We have two subcases according to the value of $\hat\chi(y_1)$.

{\bf Subcase 2.1.} $\hat\chi(y_1)\neq 0$.

In this case, (\ref{3.20}) implies that at least one of the $y_j$'s
for $2\le j\le k$ must have irrational mean index. Hence the theorem holds.

{\bf Subcase 2.2.} $\hat\chi(y_1)=0$.

Note that by Theorems 8.1.4 and 8.1.7 of \cite{Lon4} and our above
Proposition 2.5, we have $K(y_1)=1$. Since $\nu(y_1)\le 3$, it follows
from Proposition 2.6 and (\ref{2.14}):
\be 0=\hat\chi(y_1)=(-1)^{i(y_1)}(k_0(y_1)-k_1(y_1)+k_2(y_1)). \lb{3.21}\ee
By (iv) of Proposition 2.6, at most one of $k_l(y_1)$ for $l=0,\,1,\,2$
can be nonzero. Then (\ref{3.21}) yields $k_l(y_1)=0$ for $l=0,\,1,\,2$. Hence
it follows from Proposition 2.3 and Definition 2.4 that
\be C_{S^1,\; q}(\Psi_a, \;S^1\cdot u_1^p)=0,\qquad \forall p\in\N,\; q\in\Z,
   \lb{3.22}\ee
where we denote by $u_1$ the critical point of $\Psi_a$ corresponding to
$(\tau_1,\, y_1)$. In other words, $u_1^m$ is homologically invisible for
all $m\in\N$.

By Propositions 3.5 and 3.6, we can replace the term {\it infinite variationally
visible } in Definition 1.4 of \cite{LoZ1} (cf. Definition 15.3.3 of \cite{Lon4})
by {\it homologically visible}, and it is easy to check that all the results in
\cite{LoZ1} remain true under this change. Hence by Theorem 1.3 of \cite{LoZ1}
(cf. Theorem 15.5.2 of \cite{Lon4}), at least one of the $y_j$'s for
$2\le j\le k$ must have irrational mean index, i.e., we can forget $y_1$ and
consider only $y_j$'s for $2\le j\le k$, then apply that theorem. This proves
our theorem.

{\bf Case 3.} $Q_1=N_1(-1,\, 1)$.

In this case, by Theorems 8.1.4, 8.1.5 and 8.1.7 of \cite{Lon4}, we have
$$ i(y_1,\,m)=mi(y_1,\, 1)+2E\left(\frac{m\theta_1}{2\pi}\right)-2,
\quad \nu(y_1,\, m)=1+\frac{1+(-1)^m}{2},\qquad \forall m\in\N, $$
with $i(y_1,1)\in 2\Z+1$. Hence $K(y_1)=2$ by Proposition 2.5.
Because $y_1$ is non-degenerate, we have $k_l(y_1)=\dl_0^l$ for all $l\in\Z$
by (\ref{2.11}), (\ref{2.13}) and Definition 2.4. By Theorem 3.2, we have
$i(y_1)=i(y_1,1)-3\in 2\Z$ and $i(y_1^2)-i(y_1)=i(y_1,2)-i(y_1,1)\in 2\Z+1$.
Hence $k_0(y_1^2)=0$ by (v) of Proposition 2.6. Because $\nu(y_1^2)=2$, we
have $k_l(y_1^2)=0$ for $l\ge 2$. Then (\ref{2.14}) implies
$$ \hat\chi(y_1)=\frac{1+k_1(y_1^2)}{2}\neq 0. $$
Now (\ref{3.20}) implies that at least one of the $y_j$'s
for $2\le j\le k$ must have irrational mean index. Hence the theorem holds.

{\bf Case 4.} {\it $Q_1=N_1(-1,\, b)$ with $b=0,\, -1$ or $Q_1=R(\theta_2)$
with $\frac{\theta_2}{2\pi}=\frac{L}{N}\in\Q\cap(0,\,1)$ with $N>1$ and
$(L,\,N)=1$}.

Note first that if $Q_1=N_1(-1,\, b)$ with $b=0,\, -1$, then Theorems 8.1.5
and 8.1.7 of \cite{Lon4} imply that their index iteration formulae coincide
with that of a rotational matrix $R(\th)$ with $\th=\pi$. Hence in the
following we shall only consider the case $Q_1=R(\th_2)$ with
$\th_2/\pi\in (0,\,2)\cap \Q$. The same argument also shows that the theorem is
true for $Q_1=N_1(-1,-1)$.

By Theorems 8.1.4 and 8.1.7 of \cite{Lon4}, we have
\bea i(y_1,m)&=& m(i(y_1,1)-1) + 2E\left(\frac{m\theta_1}{2\pi}\right)
    + 2E\left(\frac{m\theta_2}{2\pi}\right)-3,\lb{3.23}\\
\nu(y_1,m)&=&3-2\varphi\left(\frac{m\theta_2}{2\pi}\right),
 \lb{3.24}\eea
with $i(y_1,\,1)\in 2\Z+1$ and all $m\in\N$. By Proposition 2.5, we have $K(y_1)=N$.
Note that because $y_1^m$ is non-degenerate for $1\le m\le N-1$,
$k_l(y_1^m)=\dl_0^l$ holds for $1\le m\le N-1$ by (\ref{2.11}), (\ref{2.13}) and
Definition 2.4. By Theorem 3.2, we have $i(y_1)=i(y_1,\,1)-3\in 2\Z$. Then
(\ref{2.14}) implies
\be \hat\chi(y_1)=\frac{N-1+k_0(y_1^N)-k_1(y_1^N)+k_2(y_1^N)}{N}. \lb{3.25}\ee
This follows from $\nu(y_1^m)\le 3$ for all $m\in\N$.

We have two subcases according to the value of $\hat\chi(y_1)$.

{\bf Subcase 4.1.} $\hat\chi(y_1)\neq 0$.

In this subcase, (\ref{3.20}) implies that at least one of the $y_j$'s
for $2\le j\le k$ must have irrational mean index. Hence the theorem holds.

{\bf Subcase 4.2.} $\hat\chi(y_1)=0$.

In this subcase, it follows from (\ref{3.25}) and (iv) of Proposition 2.6 that
\be k_1(y_1^N)=N-1>0. \lb{3.26}\ee

Using the common index jump theorem (Theorems 4.3 and 4.4 of
\cite{LoZ1}, Theorems 11.2.1 and 11.2.2 of \cite{Lon4}), we obtain some
$(T, m_1,\ldots,m_k)\in\N^{k+1}$ such that $\frac{m_1\theta_2}{\pi}\in\Z$
(cf. (11.2.18) of \cite{Lon4}) and the following hold by (11.2.6), (11.2.7) and
(11.2.26) of \cite{Lon4}:
\bea
i(y_j,\, 2m_j) &\ge& 2T-\frac{e(\gamma_j(\tau_j))}{2}, \lb{3.27}\\
i(y_j,\, 2m_j)+\nu(y_j,\, 2m_j) &\le& 2T+\frac{e(\gamma_j(\tau_j))}{2}-1, \lb{3.28}\\
i(y_j,\, 2m_j+1) &=& 2T+i(y_j,\,1). \lb{3.29}\\
i(y_j,\, 2m_j-1)+\nu(y_j,\, 2m_j-1)
 &=& 2T-(i(y_j,\,1)+2S^+_{\gamma_j(\tau_j)}(1)-\nu(y_j, 1)). \lb{3.30}
\eea

By  P. 340 of \cite{Lon4}, we have
\bea
&& 2S^+_{\gamma_j(\tau_j)}(1)-\nu(y_j,\,1) \nn\\
&&\qquad = 2S^+_{N_1(1,\,1)}(1)-\nu_1(N_1(1,\,1))
   +2S^+_{M_j}(1)-\nu_1(M_j)\nn\\
&&\qquad = 1 + 2S^+_{M_j}(1) - \nu_1(M_j)\nn\\
&&\qquad \ge -1, \qquad 1\le j\le k.\lb{3.31}\eea
In the last inequality, we have used the fact that the worst case for
$2S^+_{M_j}(1) - \nu_1(M_j)$ happens when $M_j=N_1(1,\,-1)^{\diamond 2}$
which gives the lower bound $-2$.

By Corollary 15.1.4 of \cite{Lon4}, we
have $i(y_j,\,1)\ge 3$ for $1\le j\le k$. Note that
$e(\gamma_j(\tau_j))\le6$ for $1\le j\le k$. Hence Theorem 10.2.4
of \cite{Lon4} yields
\bea i(y_j,\, m)+\nu(y_j,\, m)
&\le&  i(y_j, m+1)-i(y_j, 1)+\frac{e(\gamma_j(\tau_j))}{2}-1\nn\\
&\le& i(y_j, m+1)-1. \quad \forall m\in\N,\;1\le j\le k.\lb{3.32}\eea
Specially, we have
$$ i(y_j,\, m)<i(y_j,\, m+1),\qquad \forall m\in\N,\;1\le j\le k. $$
Now (\ref{3.27})-(\ref{3.30}) become
\bea
i(y_j,\, 2m_j) &\ge& 2T-3, \lb{3.33}\\
i(y_j,\, 2m_j)+\nu(y_j,\, 2m_j)-1 &\le& 2T+1,   \lb{3.34}\\
i(y_j,\, 2m_j+m) &\ge& 2T+3, \quad\forall\; m\ge 1,  \lb{3.35}\\
i(y_j,\, 2m_j-m)+\nu(y_j,\, 2m_j-m)-1 &\le& 2T-3,\quad\forall\; m\ge 1,\lb{3.36}
\eea
where $1\le j\le k$. By Proposition 2.3, we have
\bea
C_{S^1,\; q}(\Psi_a, \;S^1\cdot u_1^{2m_1})= \dl_{i(u_1^{2m_1})+1}^q\Q^{k_1(y_1^N)}
=\dl_{i(u_1^{2m_1})+1}^q\Q^{N-1},\lb{3.37}
\eea
Note that by Theorem 3.2
\be i(y_j^m)=i(y_j, m)-3,\qquad \forall m\in\N,\quad 1\le j\le k.\lb{3.38}\ee
Hence (\ref{3.23}) implies that $i(y_1^m)$ is even for all $m\in\N$.
This together with (\ref{3.35})-(\ref{3.38}) and Proposition 2.3 yield
\bea
&&C_{S^1,\; 2T-2}(\Psi_a, \;S^1\cdot u_1^{m})=0,\quad \forall m\in\N,\lb{3.39}\\
&&C_{S^1,\; 2T-4}(\Psi_a, \;S^1\cdot u_1^{m})=0,\quad \forall m\in\N,\lb{3.40}\\
&&C_{S^1,\; 2T-2}(\Psi_a, \;S^1\cdot u_j^{m})=0,\quad \forall m\neq 2m_j,\;2\le j\le k.\lb{3.41}\\
&&C_{S^1,\; 2T-4}(\Psi_a, \;S^1\cdot u_j^{m})=0,\quad \forall m\neq 2m_j,\;2\le j\le k.\lb{3.42}
\eea
In fact, by (\ref{3.35}), (\ref{3.36}) and (\ref{3.38}) for $1\le j\le k$, we have
$i(u_j^m)=i(y_j^m)\ge 2T$ for all $m>2m_j$ and
$i(u_j^m)+\nu(u_j^m)-1=i(y_j^m)+\nu(y_j^m)-1\le 2T-6$ for all $m<2m_j$.
Thus (\ref{3.41})-(\ref{3.42}) hold and (\ref{3.39})-(\ref{3.40}) hold for
$m\neq 2m_1$ by Proposition 2.3. Since $i(y_1^{2m_1})$ is even, by (\ref{3.37}),
(\ref{3.39})-(\ref{3.40}) also hold for $m = 2m_1$.

Thus by Propositions 3.5 and 3.6 we can find $p,q\in \{2,\ldots,k\}$ such that
\bea &&\Phi^\prime(u_p^{2m_p})=0,\quad \Phi(u_p^{2m_p})=c_{T-1},
\qquad C_{S^1,\; 2T-4}(\Psi_a, \;S^1\cdot u_p^{2m_p})\neq 0,\lb{3.43}\\
&&\Phi^\prime(u_q^{2m_q})=0,\quad \Phi(u_q^{2m_q})=c_{T},
\quad\qquad C_{S^1,\; 2T-2}(\Psi_a, \;S^1\cdot u_q^{2m_q})\neq 0,\lb{3.44}
\eea
where we denote also by $u_p^{2m_p}$ and $u_q^{2m_q}$ the corresponding
critical points of $\Phi$ and which will not be confused.

Note that by assumption (F) and Proposition 3.3, we have $c_{T-1}<c_{T}$.
Hence $p\neq q$ by (\ref{3.43}) and (\ref{3.44}). Then the proof of Lemma 3.1
in \cite{LoZ1}(cf. lemma 15.3.5 of \cite{Lon4}) yields
\be \hat i(y_p, 2m_p)<\hat i(y_q, 2m_q).\lb{3.45}\ee
Now if both $\hat i(y_p)\in\Q$ and $\hat i(y_q)\in\Q$ hold,
then the proof of Theorem 5.3 in \cite{LoZ1}(cf. Theorem 15.5.2 of \cite{Lon4}) yields
$$ \hat i(y_p, 2m_p)=\hat{i}(y_q, 2m_q). $$
Note that we may choose $T$ firstly such that $\frac{T}{M\hat i(y_j)}\in\N$ hold for
all $\hat i(y_j)\in\Q$ then use the proof of Theorem 5.3 in \cite{LoZ1}.
Here $M$ is the least integer in $\N$ that satisfies  $\frac{M\theta}{\pi}\in\Z$,
whenever $e^{\sqrt{-1}\theta}\in\sigma(\gamma_j(\tau_j))$ and $\frac{\theta}{\pi}\in\Q$
for some $1\le j\le k$. Hence either $\hat i(y_p)\notin\Q$ or $\hat i(y_q)\notin\Q$ holds.
This together with $\hat i(y_1)\notin\Q$ and $p, q\neq 1$ proves the theorem. \hfill\hb

{\bf Proof of Theorem 1.2.} We denote by $\{(\tau_j, y_j)\}_{1\le j\le 3}$
the three geometrically distinct
closed characteristics on $\Sg$, and by $\ga_j\equiv \gamma_{y_j}$ the associated
symplectic path of $(\tau_j,\,y_j)$ on $\Sg$ for $1\le j\le 3$. Then as in the
proof of Theorem 1.1, there exist $P_j\in \Sp(6)$ and $M_j\in \Sp(4)$ such
that
\be \ga_j(\tau_j)=P_j^{-1}(N_1(1,\,1)\dm M_j)P_j, \quad\forall\; 1\le j\le 3.
   \lb{3.46}\ee

As in P.356 of \cite{LoZ1}, if there is no $(\tau_j, y_j)$ with
$M_j=N_1(1,\,-1)^{\diamond 2}$ and $i(y_j, 1)=3$ in $\mathcal{V}_\infty(\Sigma, \alpha)$,
then $\varrho_n(\Sg)=3$. Hence we can use Theorem 1.4 of  \cite{LoZ1} (Theorem 15.5.2 of
\cite{Lon4}) to obtain the existence of at least two elliptic closed characteristics.
This proves the theorem.

It remains to show that if there exists a $(\tau_j, y_j)$  with
$M_j=N_1(1,\,-1)^{\diamond 2}$ and $i(y_j, 1)=3$ in $\mathcal{V}_\infty(\Sigma, \alpha)$,
we have at least two elliptic closed characteristics. We may assume
$M_1=N_1(1,\,-1)^{\diamond 2}$ and $i(y_1,1)=3$ without loss of generality. Note
that $(\tau_1,y_1)$ has rational mean index by Theorem 8.3.1 of {\cite{Lon4} and
Theorem 3.2.

By Theorem 1.3 of \cite{LoZ1}, we may assume that
$(\tau_2,y_2)$ has irrational mean index. Hence by Theorem 8.3.1 and Corollary 8.3.2
of \cite{Lon4}, $M_2\in \Sp(4)$ in (\ref{3.46}) can be connected to $R(\th_2)\dm Q_2$
within
$\Om^0(M_2)$ for some $\frac{\th_2}{\pi}\in\R\setminus\Q$ and $Q_2\in \Sp(2)$, where
$R(\th)=\left(\matrix{\cos\th & -\sin\th\cr
                      \sin\th &  \cos\th\cr}\right)$ for $\th\in\R$. Here we
use notations from Definition 1.8.5 and Theorem 1.8.10 of \cite{Lon4}. By
Theorem 2.7, the following identity holds
\be  \frac{\hat{\chi}(y_1)}{\hat{i}(y_1)} +\frac{\hat{\chi}(y_2)}{\hat{i}(y_2)}
+\frac{\hat{\chi}(y_3)}{\hat{i}(y_3)}=\frac{1}{2}.
\lb{3.47}\ee

Now if $Q_2$ is not hyperbolic, then both $(\tau_1, y_1)$ and $(\tau_2, y_2)$
are elliptic, so the theorem holds.

Hence it remains to consider the case that $Q_2$ is hyperbolic. Clearly
$(\tau_2, y_2)$ is non-degenerate, then it follows from (\ref{2.15}) that
$\hat\chi(y_2)\neq 0$. Hence (\ref{3.47}) implies that $\hat i(y_3)\in\R\setminus\Q$.
Now by Theorem 8.3.1 and Corollary 8.3.2 of \cite{Lon4}, $M_3\in \Sp(4)$ in
(\ref{3.46}) can be connected to $R(\th_3)\dm Q_3$ within $\Om^0(M_3)$ for some
$\frac{\th_3}{\pi}\in\R\setminus\Q$ and $Q_3\in \Sp(2)$. By the same reason as above,
it suffices to consider the case that $Q_3$ is hyperbolic.

Combining all the above, the only case we need to kick off is that
\be M_1=N_1(1,\,-1)^{\diamond 2},\quad i(y_1, 1)=3,\quad
M_2=R(\theta_2)\diamond Q_2,\quad
M_3=R(\theta_3)\diamond Q_3, \lb{3.48}\ee
where both $Q_2$ and $Q_3$ are hyperbolic.
Hence by Theorem 8.3.1 of \cite{Lon4} and Theorem 3.2, we have
\bea i(y_1^m)&=&m(i(y_1, 1)+1)-4=4m-4,\;\nu(y_1^m)=3,\quad \forall m\in\N,\lb{3.49}\\
i(y_j^m)&=&m(i(y_j)+3)+2E\left(\frac{m\theta_j}{2\pi}\right)-5,\;\nu(y_j^m)=1,
  \quad \forall m\in\N,\;j=2,3.\lb{3.50}
\eea
By Proposition 2.5, we have $K(y_1)=1$. Note that $i(y_1)=i(y_1, 1)-3=0$ by Theorem 3.2.
Hence Proposition 2.6, (\ref{2.14}) and (\ref{2.15}) imply
\bea \hat\chi (y_1)&\le& 1,\qquad \hat\chi (y_1)\in\Z,\lb{3.51}\\
\hat\chi(y_j)
    &=& \left\{\matrix{-1, & {\rm if\;\;} i(y_j)\in 2\N_0+1,  \cr
           \frac{1}{2}, &   {\rm if\;\;} i(y_j)\in 2\N_0,\cr}\right.\quad j=1,2.\lb{3.52}
\eea
By (\ref{3.49}) and (\ref{3.50}), we have
\bea \hat i(y_1)&=&4,\lb{3.53}\\
\hat i(y_j)&=&i(y_j)+3+\frac{\theta_j}{\pi}>3,\quad j=2,3.\lb{3.54}
\eea
By (\ref{3.51})-(\ref{3.54}), in order to make (\ref{3.47}) hold,
we must have
\bea \hat\chi (y_1)&=&1,\lb{3.55}\\
i(y_j)&\in&2\N_0,\quad j=2,3.\lb{3.56}
\eea
In fact, by (\ref{3.52}) and (\ref{3.54}), we have
$$\frac{\hat{\chi}(y_2)}{\hat{i}(y_2)}
+\frac{\hat{\chi}(y_3)}{\hat{i}(y_3)}<\frac{1}{6}+\frac{1}{6}<\frac{1}{2}.$$
Thus to make (\ref{3.47}) hold, we must have
$\frac{\hat{\chi}(y_1)}{\hat{i}(y_1)}>0$. Hence (\ref{3.55}) follows from (\ref{3.51}).
Now if $i(y_2)\in2\N_0+1$ or  $i(y_3)\in2\N_0+1$ holds, then by (\ref{3.52}), we have
$$ \frac{\hat{\chi}(y_1)}{\hat{i}(y_1)} +\frac{\hat{\chi}(y_2)}{\hat{i}(y_2)}
+\frac{\hat{\chi}(y_3)}{\hat{i}(y_3)}<\frac{1}{4}+\frac{1}{6}<\frac{1}{2}.
$$
Hence (\ref{3.56}) must hold.

By (\ref{2.14}), (\ref{3.49}) and (\ref{3.55}), we have
$1=\hat\chi (y_1)=k_0(y_1)-k_1(y_1)+k_2(y_1)$. Since $\nu(y_1)=3$,
by Proposition 2.6, only one of $k_0(y_1),\,k_1(y_1),\,k_2(y_1)$
can be nonzero. Hence we obtain
\be k_1(y_1)=0,\quad k_0(y_1)+k_2(y_1)=1,\lb{3.57}\ee
By Proposition 2.3, we have
\be C_{S^1,\; q}(\Psi_a, \;S^1\cdot u_j^p)=0,\quad \forall p\in\N,\;q\in2\Z+1,
   \;1\le j\le 3.\lb{3.58}\ee
In fact, by (\ref{3.49}), we have $i(y_1^m)\in2\N$ for all $m\in\N$.
Thus (\ref{3.58}) holds for $j=1$ by (\ref{2.11}), (\ref{3.57}) and Definition 2.4.
By (\ref{3.50}) and (\ref{3.56}), for $j=2, 3$, we have $i(y_j^m)\in2\N$ when $m\in2\N_0+1$ and
$i(y_j^m)\in2\N_0+1$ when $m\in2\N$. In particular, all $y_j^m$ are
non-degenerate for $m\in\N$ and $j=2, 3$. Thus (\ref{3.58}) holds for $j=2, 3$ by
(\ref{2.13}).

Note that (\ref{3.58}) implies
\be M_q=0,\quad \forall q\in2\Z+1.\lb{3.59}\ee
Together with the Morse inequality Theorem 2.8, it yields
$$ -M_{2k}-\cdots -M_2 -M_0 \ge -b_{2k}-\cdots - b_2 -b_0. $$
Thus together with the Morse inequality again, it yields
$$ b_{2k}+\cdots + b_2 + b_0\ge M_{2k}+\cdots +M_2 +M_0
   \ge b_{2k}+\cdots + b_2 + b_0, $$
for all $k\ge 0$. Therefore we obtain
\be M_q=b_q, \quad \forall q\in\Z. \lb{3.60}\ee

By (\ref{3.57}), we have two cases according to the values of $k_l(y_1)$s.

{\bf Case 1.} $k_0(y_1)=1$ and $k_2(y_1)=0$.

In this case, by Propositions 2.3, 2.5 and Definition 2.4, we have
\be \dim C_{S^1,\; q}(\Psi_a, \;S^1\cdot u_1^m)=\delta_{4m-4}^q,\quad
 \forall m\in\N,\;q\in\Z.\lb{3.61}\ee
Then by (\ref{3.60}) and (\ref{2.21}), we must have
\be C_{S^1,\; 4m-4}(\Psi_a, \;S^1\cdot u_j^p)=0,\quad \forall p,\,m\in\N,\; \;j=2,3.\lb{3.62}
\ee
By (\ref{3.60}) and (\ref{2.21}) again, $M_2=b_2=1$ implies
\be C\equiv C_{S^1,\; 2}(\Psi_a, \;S^1\cdot u_j^p)=\Q,\lb{3.63}\ee
for some $p\in\N$ and $j=2$ or $3$. If $p\ge 2$, by (\ref{3.50}), we have
\be i(y_j^p)\ge 3p+2E\left(\frac{p\theta_j}{2\pi}\right)-5\ge 3.\lb{3.64}\ee
Thus $C=0$ by Proposition 2.3. Hence $p=1$. Without loss of generality, we assume
$j=2$. Then by Proposition 2.3 and (\ref{3.63}), we have
\be i(y_2)=2.\lb{3.65}\ee
Then by (\ref{3.50}), we have
\be i(y_2^m)\ge 7,\quad\forall m\ge 2.\lb{3.66}\ee
By (\ref{3.60}) and (\ref{2.21}), $M_6=b_6=1$ implies
\be C_{S^1,\; 6}(\Psi_a, \;S^1\cdot u_j^p)=\Q,\lb{3.67}\ee
for some $p\in\N$ and $j=2$ or $3$. By (\ref{3.65}) and (\ref{3.66}),
we have $j\neq 2$, i.e., $j=3$. We must have $p=1$. In fact,
by (\ref{3.61}) and (\ref{3.63}), $y_1^m$ and $y_2^n$ already
contribute a $1$ to $M_q$ for $q=0,\,2,\,4$.
Hence by (\ref{2.21}), (\ref{3.60}) and (\ref{3.56}),
we have $i(y_3)\ge 6$, and then $i(y_3^m)\ge 15$ by (\ref{3.50}) for $m\ge 2$.
Thus $p=1$ follows from Proposition 2.3. Now we have
\be i(y_3)=6.\lb{3.68}\ee
Hence by (\ref{3.53}) and (\ref{3.55}) for $y_1$, (\ref{3.50}), (\ref{3.52}),
(\ref{3.65}) and (\ref{3.68}) for $y_2$ and $y_3$, we have
$$ \frac{\hat{\chi}(y_1)}{\hat{i}(y_1)} +\frac{\hat{\chi}(y_2)}{\hat{i}(y_2)}
+\frac{\hat{\chi}(y_3)}{\hat{i}(y_3)}
=\frac{1}{4}+\frac{1}{2(5+\frac{\theta_2}{\pi})}+\frac{1}{2(9+\frac{\theta_3}{\pi})}
<\frac{1}{2}. $$
This contradicts (\ref{3.47}) and proves Case 1.

{\bf Case 2.} $k_0(y_1)=0$ and $k_2(y_1)=1$.

The study for this case is similar to that of Case 1. Thus we are rather sketch here.

In this case, by Proposition 2.3 and Definition 2.4, we have
\be \dim C_{S^1,\; q}(\Psi_a, \;S^1\cdot u_1^m)=\delta_{4m-2}^q,\quad
    \forall m\in\N,\;q\in\Z.\lb{3.69}\ee
Then by (\ref{3.60}) and (\ref{2.21}), we must have
\bea  C_{S^1,\; 4m-2}(\Psi_a, \;S^1\cdot u_j^p)=0,\quad \forall p,\,m\in\N,\; \;j=2,3.\lb{3.70}
\eea
By (\ref{3.69}), (\ref{3.60}) and (\ref{2.21}), $M_0=b_0=1$ implies
\bea C_{S^1,\; 0}(\Psi_a, \;S^1\cdot u_j^p)=\Q,\lb{3.71}
\eea
for some $p\in\N$ and $j=2$ or $3$. By (\ref{3.64}), we have $p=1$.
Without loss of generality, we assume $j=2$. Then by Proposition 2.3
and (\ref{3.50}), we have
\be i(y_2)=0,\qquad i(y_2^m)\ge 6,\quad\forall m\ge 3.\lb{3.72}\ee
By (\ref{3.60}) and (\ref{2.21}), $M_4=b_4=1$ implies
\be C_{S^1,\; 4}(\Psi_a, \;S^1\cdot u_j^p)=\Q,\lb{3.73}\ee
for some $p\in\N$ and $j=2$ or $3$. By (\ref{3.69}) and (\ref{3.72}),
as in the verification of (\ref{3.68}), we have $j=3$ and $p=1$.
Then by Proposition 2.3, we have
\be i(y_3)=4.\lb{3.74}\ee
Hence by (\ref{3.53}) and (\ref{3.55}) for $y_1$, (\ref{3.50}), (\ref{3.52}),
(\ref{3.72}) and (\ref{3.74}) for $y_2$ and $y_3$, we have
$$ \frac{\hat{\chi}(y_1)}{\hat{i}(y_1)} +\frac{\hat{\chi}(y_2)}{\hat{i}(y_2)}
+\frac{\hat{\chi}(y_3)}{\hat{i}(y_3)}
=\frac{1}{4}+\frac{1}{2(3+\frac{\theta_2}{\pi})}+\frac{1}{2(7+\frac{\theta_3}{\pi})}
<\frac{1}{2}. $$
This contradicts (\ref{3.47}) and proves Case 2 and then the whole theorem. \hfill\hb

\medskip

\noindent {\bf Acknowledgements.} I would like to sincerely thank my
Ph. D. thesis advisor, Professor Yiming Long, for introducing me to Hamiltonian
dynamics and for his valuable help and encouragement during the writing of this
paper. I would like to say that how enjoyable it is to work with him.
I would like to sincerely thank the referee for his/her
careful reading and valuable comments and suggestions.

\bibliographystyle{abbrv}

\medskip

\end{document}